\documentclass[12pt,reqno]{amsart}
\usepackage{amssymb}
\usepackage{amsmath}
\usepackage{mathscinet}

 \rm

\renewcommand{\(}{\left(}
\renewcommand{\)}{\right)}
\renewcommand{\[}{\left[}
\renewcommand{\]}{\right]}
\newcommand{\<}{\langle}
\renewcommand{\>}{\rangle}

\newcommand{\abs}[1]{\left\lvert#1\right\rvert}
\newcommand{\norm}[1]{\left\lVert#1\right\rVert}
\newcommand{\st}{\:|\:}

\newcommand{\C}{{\mathbb{C}}}
\newcommand{\R}{{\mathbb{R}}}

\renewcommand{\phi}{\varphi}

\renewcommand{\Re}{{\mathrm{Re}}}

\renewcommand{\H}{{\mathcal{H}}}

\theoremstyle{plain}
\newtheorem{thm}{Theorem}[section]
\newtheorem{lem}[thm]{Lemma}

\theoremstyle{definition}

\theoremstyle{remark}
\newtheorem{rem}[thm]{Remark}

\title{A Generalization of a result of Vemuri}

\author{Manish Chaurasia}
\thanks{Based in part on the author's doctoral thesis (IIT (BHU)), written
under the direction of Professor M. K. Vemuri.}
 \address{Department of Mathematical Sciences, IIT (BHU), Varanasi 221005}
 \email{manishchaurasia.rs.mat18@itbhu.ac.in}
\subjclass[2010]{Primary 42C10; Secondary 42C05, 42B35, 33C45, 44A20}
\keywords{Bargmann transform, Hardy's theorem, Hermite function,
Phragm\'en-Lindel\"of principle}

\begin{document}


\begin{abstract}
Assuming that a function and its Fourier transform are dominated by
Gaussians, a sharp estimate for the rate of exponential decay of its Hermite
coefficients is obtained in terms of the variances of the dominating Gaussians.
\end{abstract}

\maketitle


\section{Introduction}\label{S:intro}


If $f \in L^1(\R)$, the {\em Fourier transform} of $f$ is defined by
\begin{equation*}
  \mathcal{F}f(\xi) = \widehat f(\xi) = \frac{1}{\sqrt{2\pi}}
  \int f(x) e^{-i\xi x} \, dx.
\end{equation*}
Let $g_a(x) = e^{-ax^2/2}$.  For $a,b>0$, let
\begin{equation*}
E(a,b) = \{ f \in L^1(\R) \st
\abs{f(x)} \le C g_a(x) \text{ and } \abs{\hat f(\xi)} \le C g_b(\xi)
\quad\text{for some $C \in \R$} \}.
\end{equation*}
Let $\phi_n$ denote the $n$-th normalized Hermite function
(see \cite[Section 2.1]{ST}). Then $\{\phi_n\}_{n=0}^\infty$ forms
an orthonormal basis for $L^2(\R)$. It follows from the Cauchy-Schwarz
inequality and Mehler's formula (see \cite[Prop 2.3]{ST}) that if for
some $t>0$
\begin{equation}\label{E:exponential-decay}
\< f, \phi_n \> = O\(e^{-2nt}\), 
\end{equation}
then $f\in E(\tanh 2rt, \tanh 2rt)$
for $0<r<1$ (see e.g. \cite{Vemuri2008hermite}). Using other methods,
Radha and Thangavelu \cite{radha2009holomorphic} proved that under the
hypothesis (\ref{E:exponential-decay}), $f$ extends to $\mathbb{C}$ as an
entire function and satisfies the estimate
$$f(x+iy) = O\( e^{-\frac{1}{2} \tanh(2rt) x^2 +
  \frac{1}{2} \coth(2rt) y^2}\),$$
for $0<r<1$, and a similar estimate holds for
$\widehat{f}$ as well. 


Conversely, what can we say about the Hermite coefficients of a
function in $E(a,b)$?

This question has been answered, but not completely. It has
connections with Hardy's uncertainty principle \cite{Hardy}
(see also \cite[Theorem 7.6]{Folland-Sitaram}),
which may be stated as follows.


\begin{thm}\label{T:hardy}
If $ab > 1$ then $E(a,b) = \{0\}$.  If $ab=1$ then $E(a,b) = \C \phi_0$.
If $ab < 1$ then $\{\phi_n\}_{n=0}^\infty \subseteq E(a,b)$.
\end{thm}


This theorem characterizes the Hermite coefficients of the
functions in $E(a,b)$ when $ab \ge 1$. However, it does not give a
precise characterization when $ab<1$. In this context, Vemuri proved the
following theorem, which sharpens the last part of the trichotomy for the
case $a=b$.
\begin{thm}\label{T:vemuri}
Let $a\in(0,1)$. If $f \in E(a,a)$ then
\begin{equation*}
\langle f,\phi_n \rangle = O\[n^{-1/4} \(\frac{1-a}{1+a}\)^{n/4}\]
\end{equation*}
for $n=1,2,\dots$, and this estimate is sharp.
\end{thm}
Shortly thereafter, Garg and Thangavelu \cite{Garg2009} generalized
Theorem \ref{T:vemuri} to the several variable case.
These results characterize Hermite coefficients of the functions in
$E(a,b)$ when $ab<1$ and $a=b$. The following theorem
is our main result, and is an extension of Theorem \ref{T:vemuri}.
\begin{thm}\label{T:estimate}
Let $a,b \in (0,\infty)$ and suppose $ab<1$.  If
$f \in E(a,b)$
then
\begin{equation*}
\< f, \phi_n \>
= O\[n^{-\frac{1}{4}}\(\frac{a+b-2ab}{a+b+2ab}\)^{n/4}\]
\end{equation*}
for $n=1,2,\dots$, and this estimate is sharp.
\end{thm}


Hardy proved Theorem \ref{T:hardy} by applying the
Phragm\'en-Lindel\"of principle to the Fourier transform of $f$.
Inspired by Hardy, Vemuri applied the Phragm\'en-Lindel\"of principle
to the Bargmann transform (see \cite[Section 2.1]{ST}) to prove
Theorem \ref{T:vemuri}. Garg and Thangavelu proved their results
using the Fourier-Wigner transform and the vector valued Bargmann
transform. Here we adapt Vemuri's technique to prove Theorem
\ref{T:estimate}.

\section{The proof}\label{hermite}

Let $\H$ denote the Hilbert space of all entire functions $F$ on $\C$
such that
\begin{equation*}
\norm{F}^2_\H = \int_{\C} \abs{F(w)}^2 \,
                        \frac{e^{-\abs{w}^2/2} \, du \, dv}{\sqrt{4\pi}}
             < \infty \quad\text{($w=u+iv$)}.
\end{equation*}
For a Schwartz class function $f$, the {\em Bargmann transform} of $f$
is defined by
\begin{equation*}
Bf(w) = \frac{e^{-w^2/4}}{\sqrt{\pi}} \int_{\R} e^{xw} e^{-x^2/2} f(x) \, dx.
\end{equation*}
%
It is shown in \cite[Section 2.1]{ST} that $B$ extends to an isometric
isomophism from $L^2(\R)$ to $\H$.

Now suppose $a,b>0$, $ab<1$ and $f\in E(a,b)$. Let $\mu=\frac{1-a}{1+a}$
and $\nu=\frac{1-b}{1+b}$. Observe that $\mu,\nu\in (-1,1)$, $\mu+\nu
\in (0,2)$, and hence $(\mu+\nu-1)^2 < 1$.
Therefore
\begin{equation*}
\mu+\nu-2\mu\nu
= \frac12\[1-(\mu+\nu-1)^2+(\mu-\nu)^2\]
> 0,
\end{equation*}
and hence
\begin{equation*}
{\frac{a+b-2ab}{a+b+2ab}}={\frac{\mu+\nu-2\mu\nu}{2-\mu-\nu}}>0.
\end{equation*}
Define
\begin{equation}\label{E:Aformulae}
A(a,b)=\sqrt{\frac{a+b-2ab}{a+b+2ab}}.
\end{equation}
\begin{lem}\label{L:main-lemma}
We have $A=A(a,b)\in (0,1)$ and there exist unique numbers
$\theta_0=\theta_0(a,b), \theta_1=\theta_1(a,b) \in (0,\pi/2)$, and
$\tau=\tau(a,b)\in (-\pi/4, \pi/4)$
with the following properties.
\begin{enumerate}
\renewcommand{\theenumi}{\alph{enumi}}
\item $\theta_0 < \tau+\pi/4 < \theta_1$.
\item $A\sin(2\theta_0-2\tau) = \mu+(1-\mu)\sin^2\theta_0$,\\
$2A\cos(2\theta_0-2\tau) = (1-\mu)\sin2\theta_0$.
\item $A\sin(2\theta_1-2\tau) = \nu+(1-\nu)\cos^2\theta_1$,\\
$2A\cos(2\theta_1-2\tau) = -(1-\nu)\sin2\theta_1$.
\end{enumerate}
\end{lem}

\begin{proof}
Clearly $A\in (0,1)$ by equation (\ref{E:Aformulae}).
Since $(\mu+\nu-2\mu\nu)(2-\mu-\nu)-(\nu-\mu)^2
= 2(1-\mu)(1-\nu)(\mu+\nu) > 0$,
$$
\frac{\nu-\mu}{\sqrt{(\mu+\nu-2\mu\nu)(2-\mu-\nu)}}\in (-1,1).
$$
Therefore there exists a unique
$\tau\in (-\pi/4,\pi/4)$ such that
\begin{equation*}
\sin 2\tau = \frac{\nu-\mu}{\sqrt{(\mu+\nu-2\mu\nu)(2-\mu-\nu)}}.
\end{equation*}
Observe that $\frac{2A\cos 2\tau}{1+\mu} > 0$ and
\begin{equation*}
\(\frac{(1-\mu)-2A\sin 2\tau}{1+\mu}\)^2+\(\frac{2A\cos 2\tau}{1+\mu}\)^2=1.
\end{equation*}
Therefore there exists a unique $\theta_0 \in (0,\pi/2)$ such that
\begin{equation*}
\cos 2\theta_0 = \frac{(1-\mu)-2A\sin 2\tau}{1+\mu}, \quad\text{and}\quad
\sin 2\theta_0 = \frac{2A\cos 2\tau}{1+\mu},
\end{equation*}
and (b) follows.

Similarly $\frac{2A\cos 2\tau}{1+\nu} > 0$ and
\begin{equation*}
\(\frac{(\nu-1)-2A\sin 2\tau}{1+\nu}\)^2+\(\frac{2A\cos 2\tau}{1+\nu}\)^2=1.
\end{equation*}
Therefore there exists a unique $\theta_1 \in (0,\pi/2)$ such that
\begin{equation*}
\cos 2\theta_1 = \frac{(\nu-1)-2A\sin 2\tau}{1+\nu}, \quad\text{and}\quad
\sin 2\theta_1 = \frac{2A\cos 2\tau}{1+\nu},
\end{equation*}
and (c) follows.

From (b), it follows that $\cos(2\theta_0-2\tau)>0$.
Since $\theta_0\in(0,\pi/2)$ and $\tau\in(-\pi/4,\pi/4)$, it follows that
$2\theta_0-2\tau \in (-\pi/2, 3\pi/2)$.  These two facts together imply that
$\theta_0-\tau < \pi/4$, whence $\theta_0 < \tau+\pi/4$.

Similarly, it follows from  (c) that
$\theta_1 > \tau+\pi/4$.  Hence (a) follows.
\end{proof}
\begin{rem}\label{R:symmetry}
Observe that $A(b,a)=A(a,b)$, $\tau(b,a)=-\tau(a,b)$,
$\theta_0(b,a)=\frac{\pi}{2}-\theta_1(a,b)$, and
$\theta_1(b,a)=\frac{\pi}{2}-\theta_0(a,b)$.
\end{rem}
\begin{lem}\label{L:Bargmann-estimate}
Let $a,b \in (0,\infty)$ and suppose $ab<1$.  Let $m=\min\{a,b\}$ and
let $A$, $\tau$, $\theta_0$
and $\theta_1$ be as in Lemma \ref{L:main-lemma}. If $f \in E(a,b)$ then
\begin{equation}\label{E:cask-strength}
\abs{Bf(w)} \le C \sqrt{\frac{2}{1+m}}
\exp\(A \frac{r^2}{4} \sin\(2\theta- 2\tau \)\), \quad (w=re^{i\theta}),
\end{equation}
for $\theta_0 \le \theta \le \theta_1$,
$\theta_0+\pi \le \theta \le \theta_1+\pi$, and
\begin{equation}\label{E:cask-strength1}
\abs{Bf(w)} \le C \sqrt{\frac{2}{1+m}}
\exp\(
A \frac{r^2}{4} \sin\(-2\theta- 2\tau \)
\), \quad (w=re^{i\theta}),
\end{equation} 
for $2\pi-\theta_1 \le \theta \le 2\pi-\theta_0$, and
$\pi-\theta_1 \le \theta \le \pi-\theta_0$.
\end{lem}
\begin{proof}
Suppose $a,b>0$, $ab<1$ and $f\in E(a,b)$. Then there exists $C>0$
such that $\abs{f(x)} \le C g_a(x)$, $x\in\R$, and $\abs{\widehat{f}(\xi)}
\le C g_b(\xi)$, $\xi\in\R$. Write $w = u+iv = re^{i\theta}$.
By arguments analogous to those used in
(\cite[Theorem 2.1]{Vemuri2008hermite}), we have
\begin{equation}\label{E:B1}
\abs{Bf(w)} \le  C \sqrt{\frac{2}{1+a}}
         	\exp \frac{(\mu + (1-\mu)\sin^2 \theta)r^2}{4},
\end{equation}
and
\begin{equation}\label{E:B2}
\abs{Bf(w)} \le  C \sqrt{\frac{2}{1+a}}
         	\exp \frac{(\nu + (1-\nu)\cos^2 \theta)r^2}{4}.
\end{equation}


Define

\begin{equation*}
F(re^{i\theta}) = \exp\( i A \frac{r^2}{4} e^{2i\(\theta-\tau\)} \) Bf(w),
\end{equation*}
where $A$ and $\tau$ are as in Lemma \ref{L:main-lemma}.
From equation (\ref{E:B1}) and Lemma \ref{L:main-lemma}(b), it follows that
\begin{equation*}
\begin{aligned}
\abs{F(re^{i\theta_0})}
\le&\; \exp\(-A \frac{r^2}{4} \sin(2\theta_0-2\tau) \)\abs{Bf(re^{i\theta_0})}\\
\le&\; C\sqrt{\frac{2}{1+a}}\exp\(\frac{r^2}{4}
        \[-A\sin(2\theta_0-2\tau)+(\mu+(1-\mu)\sin^2\theta_0)\]\)\\
\le & \; C \sqrt{\frac{2}{1+a}}.
\end{aligned}
\end{equation*}
Similarly, from equation (\ref{E:B2}) and Lemma \ref{L:main-lemma}(c),
it follows that
\begin{equation*}
\abs{F(re^{i\theta_1})} \le C \sqrt{\frac{2}{1+b}}.
\end{equation*}
Then $F$ is entire, bounded by $\sqrt{2}C e^{\abs{w}^2/2}$
everywhere (by equation (\ref{E:B1}) and the definition of $F$),
and by $C \sqrt{\frac{2}{1+m}}$ on the rays $\theta=\theta_0$ and
$\theta=\theta_1$. It follows from the Phragm\'en-Lindel\"of principle that
\begin{equation*}
\abs{F(w)} \le C \sqrt{\frac{2}{1+m}}
\end{equation*}
for $\theta_0 \le \theta \le \theta_1$
(note that $\theta_0 < \theta_1$ by Lemma \ref{L:main-lemma}(a)).
Therefore
\begin{equation*}
\abs{Bf(w)} \le C \sqrt{\frac{2}{1+m}}
\exp\(
A \frac{r^2}{4} \sin\(2\theta- 2\tau \)
\)
\end{equation*}
for $\theta_0 \le \theta \le \theta_1$.\\

Observe that $\mathcal{F}^kf\in E(a,b)$ or $E(b,a)$ according as $k$ is even
or odd.  Firstly 
\begin{equation*}
\abs{Bf(-w)}=\abs{B\mathcal{F}^2f(w)} \le C \sqrt{\frac{2}{1+m}}
\exp\(
A \frac{r^2}{4} \sin\(2\theta- 2\tau\)
\)
\end{equation*}
for $\theta_0\le\theta\le\theta_1$, or 
\begin{equation*}
\abs{Bf(w)} \le C \sqrt{\frac{2}{1+m}}
\exp\(
A \frac{r^2}{4} \sin\(2\theta- 2\tau \)
\)
\end{equation*}
for $\theta_0+\pi \le \theta \le \theta_1+\pi$.
Also
\begin{equation*}
\abs{Bf(-i w)}=\abs{B\mathcal{F}f(w)} \le C \sqrt{\frac{2}{1+m}}
\exp\(
A(b,a) \frac{r^2}{4} \sin\(2\theta- 2\tau(b,a)\)
\)
\end{equation*}
for $\theta_0(b,a) \le \theta \le \theta_1(b,a)$. By this and
Remark \ref{R:symmetry} we have
\begin{equation*}
\abs{Bf(w)} \le C \sqrt{\frac{2}{1+m}}
\exp\(
A \frac{r^2}{4} \sin\(-2\theta- 2\tau \)
\)
\end{equation*}
for $2\pi-\theta_1 \le \theta \le 2\pi-\theta_0$. Repeating the first part of
the argument with $\mathcal{F}f$ in place of $f$ shows that
(\ref{E:cask-strength1}) also holds for $\pi-\theta_1\le\theta\le \pi-\theta_0$.
\end{proof}


Now define
$\gamma_n(t) = \sqrt{\frac{2n}{A}} e^{it} \,\, \text{for} \,\, 0\le t \le 2\pi$. By the Cauchy integral formula for derivatives, we have
$Bf(w)=\sum_{n=1}^{\infty}c_n w^n$ where
$$
c_n
= \frac{1}{2\pi i} \int_{\gamma_n} \frac{Bf(w)}{w^{n+1}} \,dw.
$$
Therefore
\begin{equation}\label{E:dihedral}
\begin{aligned}
\abs{c_n}
\le&\; \frac{1}{2\pi} \int_{\gamma_n}
       \frac{\abs{Bf(w)}}{\abs{w}^{n+1}} \,\abs{dw}\\
  =&\; \frac{1}{2\pi}\(\frac{A}{2n}\)^{n/2}
     \int_{0}^{2\pi} \abs{Bf\(\sqrt{\frac{2n}{A}} e^{it}\)} \,dt\\
  =&\; \frac{1}{2\pi} \(\frac{A}{2n}\)^{n/2}
%
%
       \sum_{k=1}^4 \int_{\frac{(k-1)\pi}{2}}^{\frac{k\pi}{2}}
       \abs{Bf\(\sqrt{\frac{2n}{A}} e^{it}\)} \,dt.
\end{aligned}
\end{equation}
By Lemma \ref{L:Bargmann-estimate}, and inequalities (\ref{E:B1}) and
(\ref{E:B2}) we have
\begin{equation*}
\begin{aligned}
 \int_{0}^{\pi/2} \abs{Bf\(\sqrt{\frac{2n}{A}} e^{it}\)}\, dt
\le&\; C\sqrt{\frac{2}{1+m}} (I_n+J_n+K_n)
\end{aligned}
\end{equation*}
where
\begin{equation*}
\begin{aligned}
I_n
=&\; \int_0^{\theta_0} \exp\(\frac{(\mu+(1-\mu)\sin^2t)\, n}{2A}\) dt,\\
J_n
=&\; \int_{\theta_0}^{\theta_1} \exp\(\frac{n}{2}\sin(2t-2\tau)\) \, dt,
     \quad\text{and} \\
K_n
=&\; \int_{\theta_1}^{\frac{\pi}{2}}  \exp\(\frac{(\nu+(1-\nu)\cos^2t)\, n}{2A}\)\, dt.
\end{aligned}
\end{equation*}
By \cite[Theorem 7.7.5]{Hormander}, we have
\begin{equation*}
\begin{aligned}
J_n =&\; \int_{\theta_0}^{\theta_1}\exp\(\frac{n}{2}\)
         \exp\[i\frac{n}{2}\(i\(1-\sin(2t-2\tau)\)\)\] dt\\
=&\; O\(n^{-1/2} e^{n/2}\).
\end{aligned}
\end{equation*}
Also
\begin{equation*}
\begin{aligned}
I_n
\le&\; \theta_0\exp\(\frac{(\mu+(1-\mu)\sin^2\theta_0)\, n}{2A}\)\\
=&\; \theta_0\exp\(\frac{n}{2}\sin(2\theta_0-2\tau)\)
     \quad(\text{by Lemma \ref{L:main-lemma}(b)})\\
\le&\; \(\frac{\theta_0}{\tau+\frac{\pi}{4}-\theta_0}\)J_n. 
\end{aligned}
\end{equation*}
Similarly,
$K_n\le \(\frac{\frac{\pi}{2}-\theta_1}{\theta_1-\tau-\frac{\pi}{4}}\)J_n $
by Lemma \ref{L:main-lemma}(c).  Therefore 
$$
\int_{0}^{\pi/2} \abs{Bf\(\sqrt{\frac{2n}{A}} e^{it}\)}\, dt = O\(n^{-1/2} e^{n/2}\).
$$
The other three integrals in (\ref{E:dihedral}) are also $O\(n^{-1/2}e^{n/2}\)$
by Lemma \ref{L:Bargmann-estimate}, and the fact that the right hand sides of
inequalities (\ref{E:B1}) and (\ref{E:B2}) do not change when we replace
$\theta$ by $\pi-\theta$ or $2\pi-\theta$.
We conclude from equation (\ref{E:dihedral}) that
\begin{equation*}
c_n=O\[n^{-1/2}\(\frac{Ae}{2n}\)^{n/2}\].
\end{equation*}
Since $\abs{\<f, \phi_n\>} = \sqrt{2^n n!\pi^{1/2}} \abs{c_n}$
(see \cite{Vemuri2008hermite}),
it follows from Stirling's formula and equation
(\ref{E:Aformulae}) that
\begin{equation*}
\<f, \phi_n\> = O\(n^{-\frac{1}{4}}A^{n/2}\)
              = O\[n^{-\frac{1}{4}}\(\frac{a+b-2ab}{a+b+2ab}\)^{n/4}\].
\end{equation*}

Now we will show that this estimate is sharp.
Let
\begin{equation*}
f(x)= \exp\[-\(\frac{1+iAe^{-2i\tau}}{1-iAe^{-2i\tau}}\)\frac{x^2}{2}\]
\end{equation*}
where $A$ and $\tau$ are as in Lemma \ref{L:main-lemma}.
We claim that $f\in E(a,b)$ and
\begin{equation*}
\abs{\<f, \phi_n\>} \sim
\(\frac{2}{\pi^3}\)^{1/4} n^{-1/4} \(\frac{a+b-2ab}{a+b+2ab}\)^{n/4},
\quad n=0,2,4,\dots.
\end{equation*}

Indeed,
\begin{equation*}
\begin{aligned}
\abs{f(x)}
=&\; \exp\[-\Re \(\frac{1+iAe^{-2i\tau}}{1-iAe^{-2i\tau}}\)\frac{x^2}{2}\]
     = e^{-\frac{ax^2}{2}}, \qquad\text{and}\\
\abs{\hat{f}(\xi)}
=&\; \sqrt{2\pi}\(\frac{b}{a}\)^{1/4}
     \exp\[-\Re \(\frac{1-iAe^{-2i\tau}}{1+iAe^{-2i\tau}}\)\frac{\xi^2}{2}\]\\
=&\; \sqrt{2\pi}\(\frac{b}{a}\)^{1/4} e^{-\frac{b\xi^{2}}{2}}.
\end{aligned}
\end{equation*}
It follows that $f\in E(a,b)$.  However,
\begin{equation*}
\begin{aligned}
Bf(w)
=&\; \frac{e^{-w^2/4}}{\sqrt{\pi}}  \int e^{xw} e^{-x^2/2}
     \exp\[-\(\frac{1+iAe^{-2i\tau}}{1-iAe^{-2i\tau}}\)\frac{x^2}{2}\] \, dx\\
= & \; \frac{1}{\sqrt{\pi}} \exp\(-iAe^{-2i\tau}\frac{w^2}{4}\)\\
= & \; \frac{1}{\sqrt{\pi}} \sum_{n=0}^{\infty} \frac{(-i)^n A^n e^{-2in\tau}w^{2n}}{4^n n!}.
\end{aligned}
\end{equation*}
Therefore $\langle f, \phi_n\rangle=0$ for $n$ odd and
\begin{equation*}
\abs{\<f, \phi_n\>} = \frac{2^{-n/2}A^{n/2}\sqrt{n!}}{\sqrt{\pi}(\frac{n}{2})!},
\quad n=0,2,4,\dots.
\end{equation*}
Therefore, by Stirling's formula, we have
\begin{equation*}
\abs{\<f, \phi_n\>} \sim \(\frac{2}{\pi^3}\)^{1/4} n^{-1/4}
                         \(\frac{a+b-2ab}{a+b+2ab}\)^{n/4},
\quad n=0,2,4,\dots.
\end{equation*}

\bibliographystyle{amsplain}
\bibliography{v23-hcfahc}

\providecommand{\bysame}{\leavevmode\hbox to3em{\hrulefill}\thinspace}
\providecommand{\MR}{\relax\ifhmode\unskip\space\fi MR }
\providecommand{\MRhref}[2]{%
  \href{http://www.ams.org/mathscinet-getitem?mr=#1}{#2}
}
\providecommand{\href}[2]{#2}
\begin{thebibliography}{1}

\bibitem{Folland-Sitaram}
G.~B. Folland and A.~Sitaram, \emph{{The Uncertainty Principle: A Mathematical
  Survey}}, The Journal of Fourier Analysis and Applications \textbf{3} (1997),
  no.~3, 207--238.

\bibitem{Garg2009}
R.~Garg and S.~Thangavelu, \emph{On the {H}ermite expansions of functions from
  the {H}ardy class}, Studia Math. \textbf{198} (2010), no.~2, 177--195.
  \MR{2640076}

\bibitem{Hardy}
G.~H. Hardy, \emph{{A Theorem Concerning Fourier Transforms}}, Journal of the
  London Mathematical Society \textbf{s1-8} (1933), no.~3, 227--231.

\bibitem{Hormander}
L~H{\"o}rmander, \emph{{The analysis of linear partial differential operators.
  I}}, second ed., Grundlehren der mathematischen Wissenschaften [Fundamental
  Principles of Mathematical Sciences], vol. 256, Springer-Verlag, Berlin,
  1990, Distribution theory and Fourier analysis.

\bibitem{radha2009holomorphic}
R.~Radha and S.~Thangavelu, \emph{Holomorphic {S}obolev spaces, {H}ermite and
  special {H}ermite semigroups and a {P}aley-{W}iener theorem for the windowed
  {F}ourier transform}, J. Math. Anal. Appl. \textbf{354} (2009), no.~2,
  564--574. \MR{2515237}

\bibitem{ST}
S.~Thangavelu, \emph{Hermite and {L}aguerre semigroups: some recent
  developments}, Orthogonal families and semigroups in analysis and
  probability, S\'{e}min. Congr., vol.~25, Soc. Math. France, Paris, 2012,
  pp.~251--284.

\bibitem{Vemuri2008hermite}
M.~K. Vemuri, \emph{Hermite expansions and hardy's theorem}, arXiv preprint
  arXiv:0801.2234 (2008).

\end{thebibliography}

\end{document}